\documentclass[11pt,leqno]{article}
\usepackage{amsmath,amssymb,amsfonts}
\usepackage{graphicx}

\newtheorem{theorem}{{\sc Theorem}}
\newcommand{\bt}{\begin{theorem}}
\newcommand{\et}{\end{theorem}}
\newcommand{\newsection}[1]{\setcounter{equation}{0} \setcounter{theorem}{0}
\section{#1}}

\newcommand{\NI}{\noindent}
\newcommand{\bea}{\begin{eqnarray}}
\newcommand{\eea}{\end{eqnarray}}

\def \spec#1 {\mathop{#1}}

\def \b #1 {\bf #1}
\newcommand {\nnb }{\nonumber}

\newcommand {\CC}{\centerline}

\newcommand{\clf}{{\cal F}}

\newcommand{\cll}{{\cal L}}

\newcommand{\ity}{\infty}
\newcommand{\raro}{\rightarrow}

\newcommand{\vsp}{\vskip 1em}

\newcommand{\be}{\begin{equation}}
\newcommand{\ee}{\end{equation}}
\newcommand{\ben}{\begin{eqnarray*}}
\newcommand{\een}{\end{eqnarray*}}

\begin{document}
\CC{\bf PARAMETRIC ESTIMATION  FOR STOCHASTIC WAVE EQUATIONS}
\CC{\bf DRIVEN BY AN INFINITE DIMENSIONAL} 
\CC{\bf FRACTIONAL BROWNIAN MOTION}
\vsp
\CC{\bf B.L.S. Prakasa Rao \footnote {CR Rao Advanced Institute of Mathematics, Statistics and Computer Science, Hyderabad, India, blsprao@gmail.com}}
\vsp
\NI{\bf Abstract :} We study the problem of estimation of the parameter in a stochastic wave equation driven by an infinite dimensional fractional Brownian motion.  
\vsp
\NI{\it Keywords and Phrases:} Stochastic wave equation; Parametric estimation; fractional Brownian motion; Cylindrical fractional Brownian motion.
\vsp
\newsection{Introduction} Stochastic wave equation is used to model wave propagation through medium that exhibits random fluctuations, impurities and unpredictable external forces. Examples of such phenomena include seismic waves traveling through earth's crust which contain highly irregular rock densities and unpredictable features, acoustics to simulate sound wave transmission through the ocean accounting for random temperature gradients, salinity fluctuations and turbulent currents and for atmospheric optics to track laser beam and light propagation through atmospheric turbulence which causes random variation in the refractive index. Other examples of application include quantum field theory, thermal vibrations, ocean engineering, bridge dynamics and fiber optic networks. Statistical inference for stochastic wave equation driven by space-time white noise has been studied by Liu and Lototsky (2010) and Delgado-Vences and Pavon-Espanol (2023). Our interest is to investigate estimation of the parameters of a stochastic wave equation when it is driven by a space-time fractional Gaussian noise and study their asymptotic properties. We use the spectral approach which consists of obtaining the solution of the wave equation as a Fourier series and then using the first N modes of the solution in the estimation on a fixed window $[0,T]$. The asymptotic properties of the estimator are investigated as $N$ and $T$ tend to infinity. For reviews of earlier work on statistical inference for SPDEs, see Prakasa Rao (2001, 2010, 2026), Cialenco (2018) and Lototsky and Rozovsky (2017).
\vsp
Our aim in this paper is to initiate the study of parametric inference for stochastic wave equation driven by a cylindrical   fractional Brownian motion (fBm). The stochastic partial differential equation (SPDE) is given by 
\be
\frac{\partial^2 u}{\partial t^2}=\theta_1 \frac{\partial^2 u}{\partial x^2}+\theta_2 \frac{\partial u}{\partial t}+ \dot{W}^{H}(t), 0<t<T, 0<x<\pi
\ee
with zero initial boundary conditions where $\theta_1 \geq 1, |\theta_2|\leq 1$ driven by cylindrical fractional Brownian motion $W^H$ on $L_2[0,1]$ with the Hurst index $H\in (1/2,1)$. The solution of this equation can be written as a Fourier Series
$$u(t,x)=\sqrt{\frac{2}{\pi}}\sum_{k=1}^\ity u_k(t) \sin{kx}.$$
\vsp
We will prove the existence and uniqueness of the solution $u(t,x)$ of the equation in the next section. Our aim is to construct and study the properties of the maximum likelihood estimators of the unknown parameters $\theta_1$ and $\theta_2$ given $\{u_1(t),\dots,u_N(t), t\in [0,T]\}$ the first $N$ Fourier Coefficients of the solution $u(t,x).$

\newsection{Existence and uniqueness of solution for stochastic wave equations driven by cylindrical fBm}
\vsp 
We now extend the works of Liu and Lototsky (2010) and Delgado-Vences and Pavon-Espanol (2023) dealing with problems of estimation for stochastic wave equation  driven by additive space-time Gaussian white noise to the Stochastic Wave equation driven by additive space-time fractional Gaussian noise. We introduce some notation and  introduce the probability structure following the framework developed in Liu and Lototsky (2010). We now introduce some notation.
\vsp

Let $ \{\Omega,\clf, \clf_t, t \geq 0, P\}$ be a stochastic basis under the usual assumptions of completeness of ${\cal F}_0$ and the right-continuity of $\clf_t$. We will also suppose that  the stochastic basis supports countably many independent standard fractional Brownian motions. Recall that an infinite dimensional fractional Brownian motion $W=W^H(t), t\geq 0$ on a Hilbert space ${\cal H}$ is a linear mapping $W:f\raro W_f(.)$ from ${\cal H}$ to the space of zero-mean fractional Brownian motions such that, for every $f,g \in {\cal H}$ and $t,s >0,$ 
\vsp
$$E(W^H_f(t)W^H_g(s))= R(t,s) (f,g)_{\cal H}$$
where 
$$R(t,s)= \frac{1}{2}(|t|^{2H}+|s|^{2H}-|t-s|^{2H})$$
and $(f,g)_{\cal H}$ is the inner product of $f$ and $g$ in ${\cal H}$ and $H$ is the Hurst index of the fractional Brownian motion. If $\{e_k, k \geq 1\}$ is an orthonormal basis in ${\cal H}$ and $\{W_k^H, k \geq 1\}$ are independent standard fractional Brownian motions, then the map 
$$f \raro \sum_{k\geq 1}W^H(t) (f,e_k)_{\cal H}$$
is called a {\it cylindrical fractional Brownian motion} or an {\it infinite dimensional fractional Brownian motion} with Hurst index $H.$
A cylindrical fractional Brownian motion is represented by the generalized Fourier Series
$$W^H(t)=\sum_{k\geq 1}W_k^H(t)e_k$$
where $\{W_k(t), t\geq 0\}, k \geq 1$ are independent fractional Brownian motions.
\vsp
For $\gamma \in R,$ define the Hilbert space ${\cal H}^\gamma$ to be the closure of the set of smooth compactly supported  functions on $(0,\pi)$ with respect to the norm
$$ ||f||_\gamma= (\sum_{k\geq 1}k^{2\gamma}f_k^2)^{1/2}$$
where
$$f_k= \sqrt{\frac{2}{\pi}}\int_0^\pi f(x) \sin(kx)dx.$$
Observe that each of the functions $\sin(kx)$ belongs to every ${\cal H}^\gamma$  and if $f$ is twice continuously differentiable function with derivative $f^{\prime\prime}(x)$ on $(0,\pi)$ with $f(0)=f(\pi)=0,$ then, after integration by parts twice, it can be checked that 
$$|f_k|\leq k^{-2} \sup_{x \in (0,\pi)}|f^{\prime\prime}(x)|$$
which implies that $f\in {\cal H}^1.$ Furthermore, every $f \in {\cal H}^\gamma$ can be identified with a sequence $\{f_k, k \geq 1\}$  of real numbers such that $\sum_{k\geq 1}k^{2\gamma}f_k^2<\ity.$ Given $\gamma>0, f\in {\cal H}^{-\gamma}$ and $g\in {\cal H}^\gamma,$ we define
$$(f,g)=\sum_{k \geq 1}f_kg_k.$$
If $f,g \in L_2((0,\pi))$, then 
$$(f,g)=\int_0^\pi f(x)g(x)dx.$$
The mapping $(.,.)$, defined above, is the duality between ${\cal H}^\gamma$ and ${\cal H}^{-\gamma}$ relative to the inner product in ${\cal H}^0= L_2((0,\pi))$ (Krein et al. (1982)).
\vsp
The equation (1.1) can be interpreted as a system of two stochastic differential equations 
\be
du=v dt, dv= (\theta_1u_{xx}+\theta_2v)dt+dW^H(t) 
\ee
with the boundary conditions 
\be
u(0,x)= \frac{\partial u(t,x)}{\partial t}|_{t=0}=0, u(t,0)=u(t,\pi)=0. 
\ee
Recall the assumption that 
\be
\theta_1 \geq 1 \;\mbox{and}\; |\theta_2|\leq 1. 
\ee
{\bf Definition 2.1:} An adapted process $U\in L_2(\Omega X (0,T) X (0,\pi))$ is called a solution of the equation (2.1) if there exists an adapted process $V$ such that \\
(i) $V \in L_2(\Omega; L_2((0,T); {\cal H}^{-1}));$ and\\
(ii)for every twice continuously differentiable function $f=f(x)$ on $(0,\pi)$ with $f(0)=f(\pi)=0$, the equations
\be
(U(t,.),f)= \int_0^t(V(t,.),f)(s)ds,
\ee
and
\be
(V(t,.),f)=\int_0^t(\theta_1(U(t,.),f^{\prime\prime})+\theta_2(V(t,.),f) (s)ds + W^H_f(t) 
\ee
hold for all $t \in [0,T]$ with probability one.
\vsp
\NI{\bf Theorem 2.1:} Under the assumptions (2.2) and (2.3), the system of equations (2.1) has a unique solution $(U,V)$ and for every $\gamma <\frac{1}{2},$
\be
U\in L_2(\Omega; L_2((0,T); {{\cal H}}^\gamma)); V \in L_2(\Omega; L_2((0,T); {{\cal H}}^{\gamma-1})).
\ee
\vsp  
\NI{\bf Proof:} Let $f(x)=\sqrt{2/\pi} \sin (kx)$ in the equation (2.4) and (2.5). Let $U_k(t)= (U(t,.),f)), V_k(t)= (V(t,.),f))$ and $W^H_k=W^H_f.$ Then
\be
U_k(t)=\int_0^tV_k(s)ds, V_k(t)= -\theta_1 k^2\int_0^tU_k(s)ds+\theta_2\int_0^t V_k(s)ds+W^H_k(t). 
\ee
By assumption on the parameters $\theta_1$ and $\theta_2$, it follows that 
\be
\theta_1 k^2> (\frac{\theta_2}{2})^2
\ee
for all $k \geq 1.$ Define 
\be
\ell_k= \sqrt{ \theta_1 k^2- (\frac{\theta_2}{2})^2}.
\ee
Applying the variation of parameters formula for a linear second order equations with constant coefficients, it follows that the solution of the system of equations (2.4) and (2.5) is given by 
\be
U_k(t)=\frac{1}{\ell_k}\int_0^te^{(\theta_2/2)(t-s)} \sin(\ell_k(t-s))dW^H_k(s)
\ee
and
\be
V_k(t)=\frac{1}{\ell_k}\int_0^te^{(\theta_2/2)(t-s)} (\ell_k \cos(\ell_k(t-s))+(\theta_2/2) \sin(\ell_k(t-s)))dW^H_k(s).
\ee
Following the  inequalities for moments of integrals with respect to a fractional Brownian motion (cf. Memin et al. (2001), Prakasa Rao (2010), pp. 11-14), it follows that, there exists a number $C= C(T,H,\theta_1,\theta_2)$ such that, for all $t \in [0,T],$ 
\be
E[U_k^2(t)]\leq C(T,H,\theta_1,\theta_2) \ell_k^{-2}\;\; \mbox{and}\;\; E[V_K^2(t)]\leq C(T,H,\theta_1,\theta_2) .
\ee
Then the Gaussian processes
\be
u(t,x)=\sqrt {\frac{2}{\pi}}\sum_{k\geq 1}U_k(t) \sin(kx), v(t,x)= \sqrt {\frac{2}{\pi}}\sum_{k\geq 1}V_k(t) \sin(kx) 
\ee
satisfy the equations (2.4)-(2.6). Uniqueness of the solution follows from the completeness of the system 
$\{\sqrt{\frac{2}{\pi}} \sin(kx), k\geq 1\}$ in $L_2((0,\pi)).$
\vsp
\NI{\bf Remark:} This result can also be derived from the general theory of stochastic hyperbolic equations as pointed out by Liu and Lototsky (2010). Our proof here is modeled on the basis of Theorem 2.2 in Liu and Lotosky (2010).
\vsp
\newsection{Estimation}
Suppose that the process $u=u(t,x), v=v(t,x)$ of the equation (2.1) is observed for all $0\leq t \leq T$ and $x \in (0,\pi).$ We will discuss the estimation of parameters $\theta_1$ and $\theta_2$ from these observations. With the abuse of notation (since the fBm $W^H$ is non-differentiable almost surely), following the methods in Liu and Lototsky (2010),  we write the equation (2.1) in the original form described in (1.1):
\be
\frac{\partial^2 u}{\partial t^2}=\theta_1 \frac{\partial^2 u}{\partial x^2} +\theta_2 \frac{\partial u}{\partial t}+ \dot{W}^{H}(t), 0 < t < T, 0
<x<\pi. 
\ee
where $\dot{W}^{H}$ is interpreted as fractional white noise. We assume that the parameters $\theta_1$ and $\theta_2$ satisfy the conditions (2.3) and Theorem 2.1 implies that the system (3.1) has a unique solution  and it has the Fourier series expansion given by (2.13). We will obtain the maximum likelihood estimators of the parameters $\theta_1$ and $\theta_2$ using the observations of the $2N$-dimensional process $\{U_k(t),V_k(t), k=1,\dots,N;t \in [0,T]\}$ given by the equations (2.10) and (2.11) and study the asymptotic properties of these estimators as $ N\raro \ity.$ Recall that
\be
U_k(t)=\int_0^tV_k(s)ds, V_k(t)= -\theta_1 k^2\int_0^t  U_k(s)ds+\theta_2\int_0^t V_k(s)ds+W^H_k(t)
\ee
from the equation (2.7). Observe that
\be
dV_k(t)= F_k(V,t)dt+ dW^H_k(t)
\ee
where $F_k(V,t)=-\theta_1 k^2\int_0^t V_k(s)ds+\theta_2 V_k(t).$ Let
\be
k_H= 2H \Gamma(\frac{3}{2}-H)\Gamma(H+\frac{1}{2}),
\ee
\be
k_H(t,s)= k_H^{-1}s^{\frac{1}{2}-H}(t-s)^{\frac{1}{2}-H}, 0<t<s,
\ee
and
\be
M_k^H(t)=\int_0^t k_H(t,s) dW^H_k(s), t \geq 0.
\ee
It is well known that the process $M_k^H$ is a Gaussian martingale with the quadratic variation $<M_k^H>_t=w_t^H=\lambda_H^{-1} t^{2-2H}$ and
\be
\lambda_H= \frac{2H \Gamma(\frac{3}{2}-H)\Gamma(H+\frac{1}{2})}{\Gamma(\frac{3}{2}-H)} 
\ee
(cf. Prakasa Rao (2010),p.25). Let $C_k(t)= F_k(V,t),$  
\bea
\;\;\;\;\\\nnb
Q_k^H(t) &= & \frac{d}{dw_t^H}\int_0^t k_H(t,s) C_k(s)ds, 0\leq t \leq T\\\nnb
&=& \frac{d}{dw_t^H}\int_0^t k_H(t,s)[-\theta_1 k^2\int_0^s V_k(u)du+\theta_2 V_k(s)]ds\\\nnb
&=& -\theta_1 k^2 \frac{d}{dw_t^H}\int_0^t k_H(t,s)(\int_0^s V_k(u)du)ds+\theta_2 \frac{d}{dw_t^H}\int_0^t k_H(t,s)V_k(s)ds\\\nnb
&=& -\theta_1 k^2 D_{1k}(t)+\theta_2 D_{2k}(t)\;\;\mbox{(say)}\;\;\\\nnb
\eea
and
\be
Z_k(t)= \int_0^t Q_k^H(s) dw_s^H + M_k^H(t), 0\leq t \leq T.
\ee
From Theorem 1.19 in Prakasa Rao (2010), it follows that the process $Z_k$ is an $({\cal F}_t)$-semimartingale and the natural filtrations of the processes $Z_k$ and $W^H_k$ coincide. Let $P_k^{Z}$ and $P_k^{W^H}$ be the probability measures generated by the processes $Z_k$ and $W^H_k$ respectively on the space $C[(0,T);R]$ of continuous real-valued functions on the interval $[0,T].$ Applying Theorem 1.20 in Prakasa  Rao (2010), it follows that the measures $P_k^{Z}$ and $P_k^{W^H}$ are absolutely continuous with respect to each other and their Radon-Nikodym derivative is given by
\be
\frac{dP_k^{Z}}{dP_k^{W^H}}= \exp(\int_0^T Q_k^H(s) dM_k^H(s)-\frac{1}{2}\int_0^T [Q_k^H(s)]^2 dw_s^H).
\ee
Since the processes $\{W_k^H, k \geq 1\}$ are independent for different $k$, so are the processes $\{V_k, k \geq 1 \}.$ The probability measure $P^{V,N}$ generated by the random vector $(V_1,\dots,V_N)$ on $C([0,T]; R^N)$ is absolutely continuous with respect to the probability measure $P^{W^H,N}$ generated by the random vector $(W^H_1,\dots,W^H_N)$ on $C([0,T]; R^N)$ and their Radon-Nikodym derivative is given by
\be
\frac{dP^{V,N}}{dP^{W^H,N}}=\exp(\sum_{k=1}^N (\int_0^T Q_k^H(s) dM_k^H(s)-\frac{1}{2}\int_0^T [Q_k^H(s)]^2 dw_s^H))
\ee
and the log-likelihood ratio is given by
\bea
\;\;\;\;\\\nnb
\lefteqn{L_N(\theta_1,\theta_2)}\\\nnb 
&=& \sum_{k=1}^N (\int_0^T Q_k^H(s)dM_k^H(s)-\frac{1}{2}\int_0^T [Q_k^H(s)]^2 dw_s^H)\\\nnb
&= &\sum_{k=1}^N (\int_0^T(-\theta_1 k^2 D_{1k}(s) + \theta_2 D_{2_k}(s))dM_k^H(s) -\frac{1}{2}\int_0^T (-\theta_1 k^2 D_{1k}(s)+\theta_2
 D_{2k}(s))^2dw_s^H).\\\nnb
\eea
\vsp
Let
\be
J_{1,N}= \sum_{k=1}^Nk^4\int_0^TD_{1k}^2(t)dw_t^H, J_{2,N}= \sum_{k=1}^N\int_0^TD_{2k}^2(t)dw_t^H;
\ee
\be
J_{12,N}=\sum_{k=1}^Nk^2\int_0^TD_{1k}(t)D_{2k}(t)dw_t^H;
\ee
\be
B_{1,N}= -\sum_{k=1}^Nk^2\int_0^TD_{1k}(t)dM_k^H(t);\mbox{and}
\ee
\be
B_{2,N}=\sum_{k=1}^N\int_0^TD_{2k}(t)dM_k^H(t).
\ee
We will now consider the problem of maximum likelihood estimation of the parameters $\theta_1$ and $\theta_2$ based on the observations
$$\{U_k(t),V_k(t),k=1,\dots, N; t \in [0,T]\}.$$
The maximum likelihood estimators $\hat \theta_{1,N}$ and $\hat \theta_{2,N}$ satisfy the equations
$$\frac{\partial Z_N(\theta_1,\theta_2)}{\partial \theta_1}|_{\theta_1=\hat \theta_{1,N},\theta_2=\hat \theta_{2,N}}=0$$
and 
$$\frac{\partial Z_N(\theta_1,\theta_2)}{\partial \theta_2}|_{\theta_1=\hat \theta_{1,N},\theta_2=\hat \theta_{2,N}}=0.$$
Solving these equations, it follows that
\be
\hat \theta_{1,N}= \frac{B_{1,N}J_{2,N}+B_{2,N}J_{12,N}}{J_{1,N}J_{2,N}-J_{12,N}^2}
\ee
and 
\be
\hat \theta_{2,N}= \frac{B_{1,N}J_{12,N}+B_{2,N}J_{1,N}}{J_{1,N}J_{2,N}-J_{12,N}^2}.
\ee
\vsp
\newsection{Special Case: $\theta_2=0$} Following Delgado-Vences and Pavon-Espanol (2023), we now consider the problem of maximum likelihood estimation of the parameter $\theta_1$  assuming that the parameter $\theta_2=0.$ We consider the SPDE
\be
\frac{\partial^2 u}{\partial t^2}=\theta_1 \frac{\partial^2 u}{\partial x^2} + \dot{W}^{H}(t), 0 < t < T, 0<x<\pi 
\ee
subject to the condition $\theta_1>0$ and boundary conditions
\be
u(0,x)= \frac{\partial u(t,x)}{\partial t}|_{t=0}=0; u(t,0)=u(t,\pi)=0.
\ee
We will now study the maximum likelihood estimation of the parameter $\theta_1$ following the ideas developed earlier in Section 3 using the martingale techniques. Define  
$f_k(x) =\sqrt{2/\pi} \sin (kx)$ and note that $U_k(t)= <U(t,.),f_k>, V_k(t)=<V(t,.),f_k)>$ and $W_k^H=W^H_{f_k}.$ Then
\be
U_k(t)= \int_0^t V_k(s)ds, V_k(t)=-\theta_1 k^2\int_0^t U_k(s)ds +W^H_k(t).
\ee
Let $\tilde M_k^H$ be the martingale constructed as described earlier with the quadratic variation $w^H$ and the log-likelihood ratio in this case is given by  
\bea
Z_{N,T}(\theta_1)&=&\sum_{k=1}^N (\int_0^T(-\theta_1 k^2 D_{1k}(s)) d\tilde M_k^H(s) -\frac{1}{2}\int_0^T (-\theta_1 k^2 D_{1k}(s))^2 dw_s^H)\\\nnb
&=& \sum_{k=1}^N (\int_0^T(-\theta_1 k^2 D_{1k}(s)) d\tilde M_k^H(s) -\frac{1}{2}\int_0^T \theta_1^2 k^4 (D_{1k}(s))^2 dw_s^H).\\\nnb
\eea
Let
\be
J_{N,T}=\sum_{k=1}^N k^4 \int_0^TD_{1k}^2(s) dw_s^H,
\ee
and
\be
B_{N,T}=-\sum_{k=1}^N k^2 \int_0^T D_{1k}(s)) d\tilde M_k^H(s).
\ee
Then the log-likelihood ratio can be written in the form
\be
Z_{N,T}(\theta_1)= \theta_1 B_{N,T}-\frac{\theta_1^2}{2}J_{N,T}
\ee
and the maximum likelihood estimator (MLE) of $\theta_1$ is given by 
\be
\hat {\theta_1}_{N,T}= \frac{B_{N,T}}{J_{N,T}}.
\ee
Let $P^{V,N}$ and $P^{W,N}$ be the probability measures generated by the process $\{(U_1(t),\dots,U_N(t)); t \in [0,T]\}$ and the process $\{(W^H_1(t),\dots,W^H_N(t)); t \in [0,T]\}$ respectively on the space $C([0,T];R)$ of continuous real-valued functions on $[0,T]$ endowed with the supremum norm. Observe that the Fisher information is given by 
\bea
I_{N,T}&=& E[(\frac{\partial}{\partial \theta_1} \log \frac{dP^{V,N}}{dP^{W,N}})^2]\\\nnb
&=&- E[(\frac{\partial^2}{\partial \theta_1^2} \log \frac{dP^{V,N}}{dP^{W,N}}]\\\nnb
&=& E[J_{N,T}]\\\nnb
&=&\sum_{k=1}^N k^4\int_0^TE[D_{1k}^2(t)]dw^H_t.\\\nnb
\eea
\vsp
We will now obtain the asymptotic properties of the estimator $\hat \theta_{1,n,T}$ as $N \raro \ity$ and $T\raro \ity.$
\vsp
Let
\be
\psi_{N,T}= \sum_{k=1}^Nk^2\int_0^TD_{1k}(t)d\tilde M^H_k(t).
\ee
Observe that
\bea 
\hat \theta_{1,N,T}-\theta_1 &=&  -\frac{\psi_{N,T}}{J_{N,T}}\\\nnb
&=& -\frac{\sum_{k=1}^N k^2\int_0^T D_{1k}(t) d\tilde M^H_k(t)}{\sum_{k=1}^N k^4 \int_0^TD_{1k}^2(t)dw^H(t)}\\\nnb
\eea
\noindent{\bf Theorem 3.1:} The maximum likelihood estimator $\hat \theta_{1,N,T}$ is strongly consistent , that is,
$$\hat \theta_{1,N,T} \raro \theta_1 \;\;\mbox{a.s. as}\;\; N,T \raro \ity$$
provided
\be
\sum_{k=1}^N k^4 \int_0^TD_{1k}^2(t)dw^H(t) \raro \ity \;\;\mbox{a.s. as} \;\;T \raro \ity.  
\ee
\vsp
\NI{\bf Proof:} This result follows by observing that the process 
$$R_{N,T}= \sum_{k=1}^N k^2\int_0^T D_{1k}(t) d\tilde M^H_k(t), T\geq 0$$
is a local martingale with the quadratic variation process 
$$<R>_{N,T}=\sum_{k=1}^N k^4 \int_0^TD_{1k}^2(t)dw^H(t)$$
and applying the strong law of large numbers (cf. Liptser (1980), Prakasa Rao (1999), p.61) under the condition (4.12), it follows that 
$$\hat \theta_{1,N,T} \raro \theta_1 \;\;\mbox{a.s. as}\;\; N,T \raro \ity.$$
\vsp
\NI{\bf Theorem 3.2:} Suppose there exists a norming function $I_{N,T}$ such that
$$I_{N,T}^2<R>_{N,T}= I_{N,T}^2 \sum_{k=1}^N k^4 \int_0^T D_{1k}^2(t)dw^H(t) \raro \eta^2 \;\;\mbox {in probability as}\;\;T \raro \ity$$
where $I_{N,T} \raro 0$ almost surely as $T\raro \ity$ and $\eta$ is a random variable such that $P(\eta >0)=1.$ Then
\be
(I_{N,T}R_{N,T}, I_{N,T}^2<R>_{N,T}) \stackrel \cll \raro (\eta Z, \eta^2)\;\;\mbox{as}\; T \raro \ity
\ee
where the random variable $Z$ has the standard normal distribution and the random variables $Z$ and $\eta$ are independent.
\vsp
\NI{\bf Proof:} This result follows as a consequence of the central limit theorem for local martingales (cf. Prakasa Rao (1999), Theorem 1.49, Remark 1.47). Observe that
\be
I_{N,T}^{-1}(\hat \theta_{1,N,T}-\theta_1)= \frac{I_{N,T}R_{N,T}}{I_{N,T}^2<R>_{N,T}}.
\ee
Applying Theorem 3.2, we obtain the following result.
\vsp
\NI{\bf Theorem 3.3:} Suppose the conditions stated in the Theorem 3.2 hold. Then
\be
I_{N,T}^{-1}(\hat \theta_{1,N,T}-\theta_1)\stackrel \cll \raro  \frac{Z}{\eta} \;\; \mbox{as}\; T \raro \ity
\ee
where the random variable $Z$ has the standard normal distribution and the random variables $Z$ and $\eta$ are independent.

\noindent{\bf Remarks 3.1:} If the random variable $\eta$ is a constant with probability one, then the limiting distribution of the MLE in Theorem 3.3 is the Gaussian  distribution with mean zero and variance $\eta^{-2}.$ Otherwise, the limiting distribution is a mixture of the Gaussian distribution with mean zero and variance $\eta^{-2}$ with the mixing distribution  as that of $\eta.$
\vsp
\noindent{\bf Remarks 3.2:}The proofs of Theorems 3.1, 3.2 and 3.3 are akin to those in Prakasa Rao (2010). These results are proved by using results on limit theorems for stochastic integrals with respect to martingales (cf. Prakasa Rao (1999)). Our contribution is to reduce the problem of study in the case of fractional noise to the problem with martingale noise using the fundamental martingale as described in Prakasa Rao (2010) in the study of statistical inference for fractional diffusion processes.
\vsp
\NI{\bf Acknowledgment :} This work was supported under the scheme ``INSA Honorary Scientist" by the Indian National Science Academy (INSA) at the CR Rao Advanced Institute of Mathematics, Statistics and Computer Science, Hyderabad, India.\\
\vsp
\NI{\bf References :}
\begin{description}
\item Cialenko, I. (2018) Statistica inference for SPDEs: an overniew, {\it  Stat. Infer Stoch. Proc.}, {\bf 21}, 309-329.
\item Delgado-Vences, F. and Pavon-Espanol, J.J. (2023) Statistical inference for a stochastic equation with Malliavin calculus, {\it Stoch. Anal. Appl.}, {\bf 41}, 447-473.
\item Krein, S.G., Petunin, Yu.I., and Semenov, E.M. (1982) {\it Interpolation of linear operators}, Volume 54, {\it Translations of Mathematical Monographs}, American Mathematical Society, Providence, Rhode Island.
\item Liptser, R. (1980) A strong law of large numbers, {\it Stochastics}, {\bf 3}, 217-228.
\item Liu, W. and Lototsky, S.V. (2010) Estimating speed and damping in the stochastic wave equation, In {\it Stochastic Partial Differential Equations and Applications}, Ed. G. Da Prato and L. Tubaro, Quaderni di Matematica, Vol. 25, Dipartimento di Matematica, Seconda Universita di Napoli, pp. 191-206.
\item Lototsky, S.V. and Rozovsky, .L. (2017) {\it Stochastic Partial Differential Equations}, Springer, Heidelerg.
\item Memin, J., Mishura, Y., and Valkeila, E. (2001) Inequalities for the moments of Wiener integrals with respect to a fractional Brownian motion, {\it Stat. Probab. Lett.}, {\bf 51}, 197-206.
\item Prakasa Rao, B.L.S. (1999) {\it Semimartingales and their Statistical Inference}, CRC Press, Boca Raton and Chapman and Hall, London.
\item Prakasa Rao, B.L.S. (2001) Statistical inference  for stochastic partial differential equations, in {\it  Selected Proceedings of the Symposium on Inference for Stochastic Processes}, Ed. I.V. Basawa, C.C. Heyde and R.L. Taylor, IMS Monograph Series, Vol. 37, Institute of Mathematical Statistics, Beachwood, OH., pp. 47-70.   
\item Prakasa Rao, B.L.S. (2010) {\it Statistical Inference for Fractional Diffusion Processes}, Wiley, London.
\item Prakasa Rao, B.L.S. (2026) {\it  Advances in Statistical Inference for Processes Driven by Fractional Processes : Inference for fractional Processes}, World Scientific, Singapore.
\item Shiryayev, A.N. (1996) {\it Probability Theory}, Springer, Heidelberg.
\end{description}
\end{document}